\input amstex
\magnification=\magstep1 
\baselineskip=13pt
\documentstyle{amsppt}
\vsize=8.7truein \CenteredTagsOnSplits \NoRunningHeads
\def\per{\operatorname{per}}

\def\PER{\operatorname{PER}}

\def\dist{\operatorname{dist}}

\topmatter
\title  Computing permanents of complex diagonally dominant matrices and tensors \endtitle 
\author Alexander Barvinok  \endauthor
\address Department of Mathematics, University of Michigan, Ann Arbor,
MI 48109-1043, USA \endaddress
\email barvinok$\@$umich.edu  \endemail
\date September 2018 \enddate

\thanks  This research was partially supported by NSF Grant DMS 1361541.
\endthanks 
\keywords permanent, tensor, hypergraph, perfect matching, complex zero
\endkeywords
\abstract We prove that for any $\lambda > 1$, fixed in advance, the permanent of an $n \times n$ complex matrix, where the absolute value of each diagonal entry is at least 
$\lambda$ times bigger than the sum of the absolute values of all other entries in the same row, can be approximated within any relative error $0 < \epsilon < 1$ in quasi-polynomial
$n^{O(\ln n - \ln \epsilon)}$ time. We extend this result to multidimensional permanents of tensors and apply it to weighted counting of perfect matchings in hypergraphs.
\endabstract
\subjclass  15A15, 05C65, 41A10, 68W25, 68R05 \endsubjclass

\endtopmatter

\document

\head 1. Introduction and main results \endhead

In this paper, 
we continue the line of research started in \cite{Ba16a} and continued, in particular, in \cite{Ba17}, \cite{PR17a},  \cite{Ba16b}, \cite{PR17b}, \cite{L+17}, \cite{BR17} and \cite{EM17}, on constructing efficient algorithms for computing (approximating) combinatorially defined quantities (partition functions) by exploiting the information on their complex zeros. A typical application of the method consists of 
\smallskip
a) proving that the function in question does not have zeros in some interesting domain in ${\Bbb C}^n$ 

and

b) constructing a low-degree polynomial approximation for the logarithm of the function in a slightly smaller domain. 
\smallskip
Usually, part a) is where the main work is done: since there is no general method to establish that a multivariate polynomial (typically having many monomials) is non-zero in a domain in ${\Bbb C}^n$, some quite clever arguments are being sought and found, cf. \cite{PR17b}, \cite{EM17}, see also Section 2.5 of \cite{Ba16b} for the very few general results in this respect. Once part a) is accomplished, part b) produces a quasi-polynomial approximation algorithm in a quite straightforward way, see Section 2.2 of \cite{Ba16b}. However, if one wants to improve the complexity from quasi-polynomial to genuine polynomial time, a considerable effort can be required, see \cite{PR17a} and \cite{L+17}.

In this paper, we contribute a new method to accomplish part a) and demonstrate it by producing a quasi-polynomial algorithm to approximate permanents of complex matrices and tensors from a reasonably wide and interesting class (the hard work of sharpening our algorithms to genuine polynomial time under further restrictions on 
the matrices and tensors is done in \cite{PR17a}, see also \cite{BR17}). While computing permanents of complex matrices is of interest to quantum computations and boson sampling, see \cite{EM17}, our main contribution is to the computation of multi-dimensional permanents of (complex) tensors, which results in an efficient algorithm to count perfect matchings, weighted by their Hamming distance to one given perfect matching, in an arbitrary hypergraph. In general, the algorithm is quasi-polynomial, it becomes genuinely polynomial on hypergraphs with degrees of the vertices bounded in advance. We discuss this in Section 1.7.
 We hope that the method of this paper will find other applications, in particular to count solutions, weighted by their Hamming distance to one given solution, of other NP-complete problems.

We start with a popular example of the 
permanent of a square matrix.
\subhead (1.1) Permanents \endsubhead  
Let $A=\left(a_{ij}\right)$ be $n \times n$ complex matrix and let
$$\per A =\sum_{\sigma \in S_n} \prod_{i=1}^n a_{i \sigma(i)}$$
be its {\it permanent}.  Here $S_n$ is the symmetric group of all $n!$ permutations $\sigma$ of the set $\{1, \ldots, n\}$. First, we prove the following result.
\proclaim{(1.2) Theorem} Let $A$ be $n \times n$ complex matrix such that 
$$\sum_{j=1}^n |a_{ij}| \ < \ 1 \quad \text{for} \quad i=1, \ldots, n.$$
Then 
$$\per (I +A) \ne 0,$$
where $I$ is the $n \times n$ identity matrix.
\endproclaim
In \cite{Br59}, Brenner obtains a family of inequalities satisfied by (determinantal) minors of a diagonally dominant matrix, mentions as a corollary (Corollary 5 from \cite{Br59}), that the determinant of such a matrix is necessarily non-zero, and concludes the paper with the following sentence: ``The referee remarked that since the permanent of matrix can be expanded by minors, corresponding theorems holds for permanents."
The permanental version of Corollary 5 of \cite{Br59} is equivalent to our Theorem 1.2. 

In this paper, we provide a different proof of Theorem 1.2, which easily extends to multi-dimensional permanents of tensors (Theorem 1.5 below) and, more, generally, can be useful for establishing a zero-free region for an arbitrary multi-affine polynomial.

We also note that $A$ is Hermitian then $I+A$, being a diagonally dominant matrix, is necessarily positive definite and hence $\per (I+A)$ is positive real \cite{MN62}.

Let $A$ be a matrix satisfying the conditions of Theorem 1.2. Then we can choose a continuous branch of the function 
$A \longmapsto \ln \per(I +A)$. Applying the methods developed in \cite{Ba16b} and \cite{Ba17}, we obtain the following result.

\proclaim{(1.3) Theorem} Let us fix a real $0 < \lambda < 1$. Then for any positive integer $n$ and for any $0 < \epsilon < 1$ there exists a polynomial 
$p(A)=p_{\lambda, n, \epsilon}(A)$ in the entries of an $n \times n$ complex matrix $A$ such that $\deg p =O_{\lambda}(\ln n -\ln \epsilon)$ and 
$$\left| \ln \per(I+A) - p(A) \right| \ \leq \ \epsilon$$
provided $A=\left(a_{ij}\right)$ is a complex $n \times n$ matrix such that 
$$\sum_{j=1}^n \left| a_{ij}\right| \ < \ \lambda \quad \text{for} \quad i=1, \ldots, n.$$
\endproclaim
\noindent Moreover, given $\lambda, n$ and $\epsilon$, the polynomial $p_{\lambda, n, \epsilon}$ can be computed in $n^{O_{\lambda}(\ln n - \ln \epsilon)}$
(the implied constant in the ``$O$" notation depends only on $\lambda$).

We note that the value of $|\per(I +A)|$ for a matrix $A$ satisfying the conditions of Theorem 1.3  can vary in an exponentially wide range in $n$, even when $A$ is required to have zero diagonal: choosing $A$ to be 
block-diagonal with $2 \times 2$ blocks
$$\left( \matrix 0 & a \\ b & 0 \endmatrix \right),$$
 we can make $|\per(I+A)|$ as large as $(1+\lambda^2)^{\lfloor n/2 \rfloor}$ and as small as $(1-\lambda^2)^{\lfloor n/2 \rfloor}$.

If $B=\left(b_{ij}\right)$ is a strongly diagonally dominant complex matrix such that 
$$\lambda |b_{ii}| \ \geq \ \sum_{j: \ j \ne i} |b_{ij}| \quad \text{for} \quad i=1, \ldots, n$$
then
$$\per B =  \left(\prod_{i=1}^n b_{ii} \right)\per C,$$
where $C=\left(c_{ij}\right)$ is obtained from $B$ by a row scaling
$$c_{ij} = b_{ii}^{-1} b_{ij} \quad \text{for all} \quad i, j.$$
We have $C=I +A$, where $A$ satisfies the conditions of Theorem 1.3. Hence our algorithm can be applied to approximate the permanents of complex strongly diagonally dominant matrices.

If the matrix $A$ in Theorem 1.3 is Hermitian, then $I+A$, being diagonally dominant, is positive definite. A polynomial time algorithm approximating permanents of 
$n \times n$ positive semidefinite matrices within a simply exponential factor of $c^n$ (with $c \approx 4.84$) is constructed in \cite{A+17}. Other interesting classes
of complex matrices where efficient permanent approximation algorithms are known are some random matrices \cite{EM17} and matrices not very far from the matrix filled with 1s \cite{Ba17}. Famously, there is a randomized polynomial time algorithm to approximate $\per A$ if $A$ is non-negative real \cite{J+04}. The best known deterministic polynomial time algorithm approximates the permanent of an $n \times n$ non-negative real matrix within an exponential factor of $2^n$ and is conjectured to approximate it within a factor of $2^{n/2}$ \cite{GS14}.

\subhead (1.4) Multidimensional permanents \endsubhead For $d \geq 2$, let $A=\left(a_{i_1 \ldots i_d}\right)$ be 
a cubical $n \times \ldots \times n$ array (tensor) of complex numbers (so that for $d=2$ we obtain an $n \times n$ square matrix). 
We define the permanent of $A$ by 
$$\PER A=\sum_{\sigma_2, \ldots, \sigma_d \in S_n} \prod_{i=1}^n a_{i \sigma_2(i) \ldots \sigma_d(i)}.$$
Clearly, our definition agrees with that of Section 1.1 for the permanent of a matrix. Just as the permanent of a matrix counts perfect matchings in the underlying weighted bipartite graph, the permanent of a tensor counts perfect matchings in the underlying weighted $d$-partite hypergraph, see, for example, Section 4.4 of \cite{Ba16b}. 

We define the {\it diagonal} of a tensor $A$ as the set of entries $\left\{ a_{i \cdots i},\ i=1, \ldots, n \right\}$.
We prove the following extension of Theorem 1.2. 
\proclaim{(1.5) Theorem} Let $A=\left(a_{i_1 \ldots i_d}\right)$ be a $d$-dimensional $n \times \ldots \times n$ complex tensor with such that 
$$\sum_{1 \leq i_2, \ldots, i_d \leq n} \left| a_{i i_2 \ldots i_d} \right| \ < \ 1 \quad \text{for} \quad i=1, \ldots, n.$$
Then
$$\PER(I + A) \ne 0,$$
where $I$ is the $d$-dimensional $n \times \ldots \times n$ tensor with diagonal entries equal to 1 and all other entries equal to 0.
\endproclaim 

Similarly to Theorem 1.3, we deduce from Theorem 1.5 the following result:

\proclaim{(1.6) Theorem} Let us fix a positive integer $d \geq 2$ and a real $0 < \lambda < 1$. Then for any positive integer $n$ and any $0 < \epsilon < 1$ there exists a polynomial 
$p(A)=p_{d, \lambda, n, \epsilon}(A)$ in the entries of a $d$-dimensional $n \times \ldots \times n$ tensor $A$ such that 
$\deg p=O_{\lambda}(\ln n - \ln \epsilon)$ and 
$$\left| \ln \PER(I+A) -p(A)\right| \ \leq \ \epsilon$$
for any $d$-dimensional $n \times \ldots \times n$ tensor $A$ such that 
$$\sum_{1 \leq i_2, \ldots, i_d \leq n} \left| a_{i i_2 \ldots i_d} \right| \ < \ \lambda \quad \text{for} \quad i=1, \ldots, n.$$
\endproclaim 
Moreover, given $d$, $\lambda$, $n$ and $\epsilon$, the polynomial $p_{d, \lambda, n, \epsilon}$ can be computed in 
$n^{O_{d, \lambda}(\ln n -\ln \epsilon)}$ time (so that the implied constant in the ``$O$" notation depends only on $d$ and $\lambda$). 

\subhead (1.7) Weighted counting of perfect matchings in hypergraphs \endsubhead
We describe an application of Theorem 1.6 to weighted counting of perfect matchings in hypergraphs, cf. \cite{BR17}.
Let $H$ be a $d$-partite hypergraph with set $V$ of $nd$ vertices equally split among $d$ pairwise disjoint parts $V_1, \ldots, V_d$ such that $V=V_1 \cup \ldots \cup V_d$. 
 The set $E$ of edges of $H$ consists of some $d$-subsets of $V$ containing exactly one vertex from each part $V_i$. We number the vertices in each part $V_i$ by $1, \ldots, n$ and encode $H$ by a $d$-dimensional tensor $A=\left(a_{i_1 \ldots i_d}\right)$, where 
$$a_{i_1 \ldots i_d} = \cases 1 &\text{if\ } \{i_1, \ldots, i_d\} \in E, \\ 0 &\text{otherwise.} \endcases$$
A {\it perfect matching} in $H$ is a collection of $n$ edges containing each vertex exactly once. As is known, for $d \geq 3$ it is an NP-complete problem to determine whether a given $d$-partite hypergraph contains a perfect matching. Suppose, however, that we are given 
{\it one} perfect matching $M_0$ in $H$. Without loss of generality we assume that $M_0$ consists of the edges 
$(1, \ldots, 1)$, $(2, \ldots, 2), \ldots, (n, \ldots, n)$. For a perfect matching $M$ in $H$, let $\dist(M, M_0)$ be the number of edges in which $M$ and $M_0$ differ (the Hamming distance between $M$ and $M_0$).  Let us choose a $\lambda >0$. It is not hard to see that 
$$\PER\left(I + \lambda^2 (A-I)\right) = \sum_M \lambda^{\dist(M, M_0)}, \tag1.7.1$$
where the sum is taken over {\it all} perfect matchings $M$ in $H$. Let us assume now that each vertex of the part $V_1$ is contained in at most 
$\Delta$ edges of $H$. It follows from Theorem 1.6 that (1.7.1) can be efficiently approximated for any $\lambda$, fixed in advance, provided 
$$\lambda < \sqrt{1 \over {\Delta-1}}.$$
This is an improvement compared to \cite{BR17} where we could only afford 
$\lambda = O(1/\Delta \sqrt{d})$ and also required the degree of {\it every} vertex of $H$ not to exceed $\Delta$. As is discussed in \cite{BR17}, see also \cite{PR17a}, for any $\Delta$, fixed in advance, we obtain a polynomial time algorithm
approximating (1.7.1).  

Generally, knowing one solution in an NP-complete problem does not help one to find out if there are other solutions. Our result shows that some statistics over the set of all solutions can still be computed efficiently. 

We prove Theorems 1.2 and 1.5 in Section 2 and deduce Theorems 1.3 and 1.6 from them in Section 3.

\head 2. Proofs of Theorems 1.2 and 1.5 \endhead

We start with a simple lemma.
\proclaim{(2.1) Lemma} Let us fix $ \alpha_1, \ldots, \alpha_n \in {\Bbb C}$. Then for any $z_1, \ldots, z_n \in {\Bbb C}$ there 
exist $z_1^{\ast}, \ldots, z_n^{\ast} \in {\Bbb C}$ such that 
$$\sum_{k=1}^n \alpha_k z_k = \sum_{k=1}^n \alpha_k z_k^{\ast}, \quad \sum_{k=1}^n \left| z_k^{\ast}\right| \ \leq \ \sum_{k=1}^n \left| z_k\right|$$
and $z_k^{\ast} \ne 0$ for at most one $k$.
\endproclaim
\demo{Proof} Without loss of generality, we assume that $\alpha_k \ne 0$ for $k=1, \ldots, n$. Given 
$z_1, \ldots, z_n$, let us define 
$$K=\left\{ \left(x_1, \ldots, x_n\right) \in {\Bbb C}^n: \quad \sum_{k=1}^n \alpha_k x_k = \sum_{k=1}^n \alpha_k z_k \quad 
\text{and} \quad \sum_{k=1}^n \left| x_k\right| \ \leq \ \sum_{k=1}^n \left| z_k \right| \right\}.$$
Then $K$ is a non-empty compact set and the continuous function 
$$\left(x_1, \ldots, x_n \right) \longmapsto \sum_{k=1}^n |x_k|$$
attains its minimum on $K$ at some point, say $\left(y_1, \ldots, y_n \right)$. We claim that all non-zero complex numbers among 
$\alpha_k y_k$ are positive real multiples of each other.

Suppose that, say $\alpha_1 y_1 \ne 0$ and $\alpha_2 y_2 \ne 0$ are not positive real multiples of each other and let $v = \alpha_1 y_1+ \alpha_2 y_2$. Hence $|\alpha_1 y_1| + |\alpha_2 y_2| > |v|$. Let 
$$y_1'={v |\alpha_1 y_1| \over \alpha_1 \left(|\alpha_1 y_1| +|\alpha_2 y_2|\right)}  \quad \text{and} \quad y_2'={v |\alpha_2 y_2| \over \alpha_2 \left(|\alpha_1 y_1| +|\alpha_2 y_2|\right)} .$$
 Then 
 $$\alpha_1 y_1' +\alpha_2 y_2'=v = \alpha_1 y_1 + \alpha_2 y_2$$ and
 $$|y_1'| + |y_2'| ={|v| |y_1| \over  |\alpha_1 y_1| + |\alpha_2 y_2|} +{|v| |y_2|  \over  |\alpha_1 y_1| + |\alpha_2 y_2|} \ < \ |y_1|+|y_2|.$$
 Hence defining $y_k'=y_k$ for $n > 2$, we obtain a point $\left(y_1', \ldots, y_n'\right) \in K$ with 
 $$\sum_{k=1}^n \left| y_k'\right| \ < \ \sum_{k=1}^n \left| y_k \right|,$$
 which is a contradiction.

This proves that there is a point $\left(y_1, \ldots, y_n\right) \in K$ where all non-zero complex numbers $\alpha_k y_k$ are positive real multiples of each other,
so that 
$$\text{the ratios} \quad {\alpha_k y_k \over \left|\alpha_k y_k\right|} \quad \text{when} \quad \alpha_k y_k \ne 0 \quad \text{are all equal}.$$
Next, we successively reduce the number of non-zero coordinates among $y_1, \ldots, y_n$, while keeping all non-zero numbers $\alpha_k y_k$ positive real multiples of each other.

 Suppose that there are two non-zero coordinates, say $y_1$ and $y_2$. Without loss of generality, we assume that $|\alpha_1| \leq |\alpha_2|$.
Now, we let:
$$y_1'=0 \quad \text{and} \quad y_2'= y_2 + { |\alpha_1| |y_1| \over |\alpha_2| |y_2|} y_2.$$
Then 
$$|y_1'| + |y_2'| =|y'_2| = |y_2| + {|\alpha_1| |y_1| \over |\alpha_2|} \ \leq \ |y_1| +|y_2|.$$
Moreover,
$$\split \alpha_1 y_1' + \alpha_2 y_2'=&\alpha_2 y_2 + \alpha_2  { |\alpha_1| |y_1| \over |\alpha_2| |y_2|} y_2=\alpha_2 y_2 + {\alpha_2 y_2 \over |\alpha_2 y_2|} | \alpha_1 y_1| =
\alpha_2 y_2 + {\alpha_1 y_1 \over |\alpha_1 y_1|} |\alpha_1 y_1|\\=&\alpha_2 y_2 + \alpha_1 y_1. \endsplit$$
Hence letting $y_k'=y_k$ for $k > 2$ we obtain a point $\left(y_1', \ldots, y_n'\right) \in K$ with fewer non-zero coordinates. Moreover, all non-zero numbers $\alpha_k y_k'$
remain positive real multiples of each other.
 Repeating this process, we 
obtain the desired vector $\left(z_1^{\ast}, \ldots, z_n^{\ast} \right) \in {\Bbb C}^n$.
{\hfill \hfill \hfill}\qed
\enddemo

Now we are ready to prove Theorem 1.2.
\subhead (2.2) Proof of Theorem 1.2 \endsubhead First, we observe that without loss of generality, we may assume that $A$ has zero diagonal. Indeed, let $A=\left(a_{ij}\right)$ be an $n \times n$ complex matrix with sums of the absolute values of entries in each row less than 1. Then the diagonal entry in the $i$-th row of the matrix $I+A$ is $1+a_{ii}$ with absolute value 
$|1+a_{ii}| \geq 1- |a_{ii}| >0$, while the sum of the absolute values of the off-diagonal entries in the $i$-th row of $I+A$ is less than $1-|a_{ii}|$. Consequently, dividing the 
$i$-row of $I+A$ by $1+a_{ii}$ for $i=1, \ldots, n$, we obtain the matrix $I + A'$ where $A'$ satisfies the conditions of the theorem and, additionally, has zero diagonal. 
Moreover,
$$\per (I + A) = \left( \prod_{i=1}^n (1+a_{ii}) \right) \per (I + A').$$
Thus we assume that $A$ has zero diagonal.

Let ${\Cal M}_n$ be the set of $n \times n$ complex matrices $A$ with zero diagonal and sums of absolute 
values of entries in every row less than 1.
We claim for every $A \in {\Cal M}_n$ there exists $B \in {\Cal M}_n$  such that $\per(I +A) =\per(I +B)$ and $B$ has at most one non-zero entry in every row.
Given $A=\left(a_{ij}\right)$, we construct the matrix $B$ step by step by modifying $A$ row by row in $n$ steps. Let $A_{ij}$ be the $(n-1) \times (n-1)$ matrix obtained from $A$ by crossing the $i$-th row and the $j$-th column.
We have 
$$\per (I +A) = \per A_{11} + \sum_{j=2}^n a_{1j} \per A_{1j}.$$
Applying Lemma 2.1, we find $b_{1j}$ for $j=1, \ldots, n$ such that $b_{11}=0$, 
$$\sum_{j=2}^n b_{1j} \per A_{1j} =\sum_{j=2}^n a_{1j} \per A_{1j}, \quad   \sum_{j=2}^n \left| b_{1j}\right|  \ < \ 1$$
and at most one of the numbers $b_{1j}$ is non-zero. At the first step, we define $B$ by replacing $a_{1j}$ by $b_{1j}$ for $j=1, \ldots, n$ and note that $\per(I+A)=\per(I+B)$.

At the end of the $(k-1)$-st step, we have a matrix $B \in {\Cal M}_n$ such that $\per(I +A)=\per(I+B)$ and each of the first $k-1$ rows of $B$ contains at most one non-zero entry.
If $k \leq n$, we write 
$$\per (I +B) = \per B_{kk}+  \sum_{j:\ j\ne k}  a_{kj} \per B_{kj},$$
where $B_{kj}$ is the $(n-1) \times (n-1)$ matrix obtained from $B$ by crossing out the $k$-th row and $j$-th column. Applying Lemma 2.1, we find $b_{kj}$ for $j=1, \ldots, n$
such that $b_{kk}=0$, 
$$\sum_{j:\ j \ne k} b_{kj} \per B_{kj} = \sum_{j:\ j \ne k} a_{kj} \per B_{kj}, \quad \sum_{j:\ j \ne k} \left|b_{kj} \right| \ < \ 1$$
and at most one of the numbers $b_{kj}$ is non-zero. We modify $B$ by replacing $a_{kj}$ with $b_{kj}$ for $j=1, \ldots, n$.
We have $\per(I + A) =\per(I+B)$.

At the end of the $n$-th step, we obtain a matrix $B \in {\Cal M}_n$ containing at most one non-zero entry in each row and such that 
$\per (I +A) =\per (I +B)$. 

We now have to prove that $\per (I+B) \ne 0$. Since every row of $B$ contains at most one non-zero entry, the total number of non-zero entries in $B$ is at most $n$.
Therefore, if there is a column of $B$ containing more than one non-zero entry, there is a column, say the $k$-th, filled by zeros only.
Then $\per (I + B) = \per(I+B')$, where $B'$ is the matrix obtained from $B$ by crossing out the $k$-th row and column. Hence without loss of generality, we may assume that every row and every column of $B$ contains exactly one non-zero entry. Let us consider the bipartite graph on $n + n$ vertices where the $i$-th vertex on one side is connected by an edge to the $j$-th vertex on the other side if and only if the $(i,j)$-th entry of $I+B$ is not zero. Then the graph is a disjoint union of some 
even cycles $C_1, \ldots, C_m$. Hence
$$\per (I+B)=\prod_{k=1}^m \left(1 + \prod\Sb \{i, j\} \in C_k \\ i \ne j \endSb  b_{ij}\right)$$
where $b_{ij}$ are the non-zero entries of $B$ corresponding to the edges in $C_k$. Since $|b_{ij}| < 1$, 
we have $\per(I+B) \ne 0$.
{\hfill \hfill \hfill} \qed

Before we prove Theorem 1.5, we introduce a convenient definition.
\subhead (2.3) Slice expansion of the permanent of a tensor \endsubhead
Let $A=\left(a_{i_1 \ldots i_d}\right)$ be a $d$-dimensional $n \times \ldots \times n$ tensor. Let us fix $1 \leq k \leq d$ and 
$1 \leq j \leq n$. We define the $(k, j)$-th {\it slice} of $A$ as the set of $n^{d-1}$ entries $a_{i_1 \ldots, i_d}$ with $i_k = j$.
In particular, if $d=2$, so $A$ is a matrix, the $(k, j)$-th slice of $A$ is the $j$-th row if $k=1$ and the $j$-th column if $k=2$.
Theorem 1.5 asserts that $\PER(I + A) \ne 0$, if $A$ is a complex tensor with sums of the absolute values of the entries in the $(1, j)$-th slice less than 1 for $j=1, \ldots, n$.

For a given entry $a_{i_1 \ldots i_d}$, let $A_{i_1 \ldots i_d}$ be the $d$-dimensional $(n-1) \times \ldots \times (n-1)$ tensor obtained from $A$ by crossing out the $d$ slices containing $a_{i_1 \ldots i_d}$. Then, for any $k=1, \ldots, d$ and any $j=1, \ldots, n$, we have the $(k, j)$-slice expansion of the permanent: 
$$\PER A = \sum\Sb i_1, \ldots, i_d: \\ i_k = j \endSb a_{i_1 \ldots i_d} \PER A_{i_1 \ldots i_d}. \tag2.3.1$$
\subhead (2.4) Proof of Theorem 1.5 \endsubhead The proof closely follows that of Theorem 1.2 in Section 2.2. 

Arguing as in Section 2.2, without loss of generality we assume that $A$ has zero diagonal.

Now we prove that for every $d$-dimensional $n \times \ldots \times n$ tensor $A$ with zero diagonal and sums of absolute values in 
the $(1, j)$-th slice less than 1 for $j=1, \ldots, n$, there exists a $d$-dimensional $n \times \ldots \times n$ tensor $B$ with zero diagonal such that $\PER(I+A)=\PER(I+B)$, all entries of $B$ are less than 1 in the absolute value and the $(1, j)$-th slice of $B$ contains at most one non-zero entry for $j=1, \ldots, n$.

We construct $B$ by modifying $A$ slice by slice starting with $B=A$. If the $(1, j)$-th slices for $j=1, \ldots, m-1$ are modified, we consider the 
$(1, m)$-slice expansion (2.3.1) 
$$\PER (I + B) = \PER B_{m \ldots m} + \sum\Sb i_2, \ldots, i_d: \\ (i_2, \ldots, i_d) \ne (m, \ldots, m) \endSb
a_{m i_2 \ldots i_d} \PER B_{mi_2 \ldots i_d}.$$
Using Lemma 2.1, we  find complex numbers $b_{m i_2 \ldots i_d}$ such that $b_{m \ldots m}=0$, at most one 
of $b_{m i_2 \ldots i_d}$ is non-zero, less than 1 in the absolute value and 
$$\sum\Sb i_2, \ldots, i_d: \\ (i_2, \ldots, i_d) \ne (m, \ldots, m) \endSb
a_{m i_2 \ldots i_d} \PER B_{mi_2 \ldots i_d} = \sum\Sb i_2, \ldots, i_d: \\ (i_2, \ldots, i_d) \ne (m, \ldots, m) \endSb
b_{m i_2 \ldots i_d} \PER B_{mi_2 \ldots i_d}.$$
We then replace the entries $a_{mi_2 \ldots i_d}$ in the $(1,m)$-th slice of $B$ by the entries $b_{mi_2 \ldots i_d}$.
In the end we produce the desired tensor $B$.

Hence our goal is to prove that $\PER(I+B) \ne 0$ if $B$ is a $d$-dimensional $n \times \ldots \times n$ tensor with zero diagonal, complex entries less than 1 in the absolute value and containing at most one non-zero entry in the $(1, j)$-th slice for 
$j =1, \ldots, n$. In particular, the total number of non-zero entries of $B$ does not exceed $n$. It follows then that if some $(k, j)$-th slice of $B$ contains more than one non-zero entry, there is a $(k, m)$-th slice of $B$ filled entirely by zeros. Denoting by $B'$ the $(n-1) \times \ldots \times (n-1)$ tensor obtained from $B$ by crossing out the $d$ slices containing $(m, \ldots, m)$, we will have $\PER(I+B) = \PER(I+B')$. Hence without loss of generality we assume that $B$ contains exactly one non-zero entry in each of the $dn$ slices. 

Let us consider the underlying $d$-partite hypergraph $H$ encoded by $I+B$. We number the vertices in each of the $d$ parts by $1, \ldots, n$. 
A set $E$ of $d$ vertices of $H$ is an edge of $H$ if and only if it contains exactly one vertex $i_k$ for $k=1, \ldots, d$ from each part and the $(i_1, \ldots, i_d)$-th entry of $I+B$ is non-zero. Hence each vertex is contained in exactly two edges of $H$. We call the edges $e_j=(j, \ldots, j)$ for $j=1, \ldots, n$ {\it standard} (they correspond to the diagonal $1$s of $I+B$) and the remaining $n$ edges {\it additional} (they correspond to the non-zero entries of $B$). We define the {\it weight} of a standard edge equal to 1 and of an additional edge 
$\{i_1, \ldots, i_d\}$ equal to $b_{i_1 \ldots i_d}$. We define the {\it weight} of a perfect matching $M$ in $H$ as the product of weights of the edges in $M$ and the total weight $w(H)$ as the sum of weights of all perfect matchings in $H$. Hence $\PER(I+B)=w(H)$ and our goal is to prove that $w(H) \ne 0$. 

The hypergraph $H$ splits into the disjoint union of {\it connected components} \newline $H_1, \ldots, H_m$, so that for any two vertices
$v'$ and $v''$ in each $H_k$ there is a chain of edges $u_1, \ldots, u_s$ such that $v' \in u_1$, $v'' \in u_s$ and $u_i \cap u_{i+1} 
\ne \emptyset$ for $i=1, \ldots, s-1$ and for any two vertices $v'$ and $v''$ in different $H_i$ and $H_j$ there is no such chain.
Then $w(H)=w(H_1) \cdots w(H_k)$ and hence it suffices to prove that $w(H_k) \ne 0$ for $k=1, \ldots, m$. Hence without loss of generality, we assume that $H$ is connected. There is a perfect matching of weight 1 in $H$ consisting of the standard edges $e_1, \ldots, e_n$, which we call {\it standard}. Furthermore, $n$ additional edges are necessarily pairwise disjoint (since otherwise there is a vertex that belongs to at least three edges: two additional and one standard) and also form a perfect matching, which we call {\it additional}.

We claim that any perfect matching in $H$ is either standard or additional. Seeking a contradiction, suppose that there is a perfect matching consisting of 
some $1 \leq k < n$ standard edges, say $e_1, \ldots, e_k$ and some $n-k > 0$ additional edges edges, say $u_{k+1}, \ldots, u_n$. Let $U=u_{k+1} \cup \ldots \cup u_n$ and let $W=e_1 \cup \ldots \cup e_k$. Then every vertex in $U$ is contained in one standard edge among $e_{k+1}, \ldots, e_n$ and one additional edge among $u_{k+1}, \ldots, u_n$  and 
hence there cannot be an edge containing a vertex of $U$ and a vertex of $W$, which contradicts the assumption that $H$ is connected. 

Hence there are exactly two perfect matchings in connected $H$, one standard of weight 1 and one additional of weight less than 1 in the absolute value. Thus 
$w(H) \ne 0$ and the proof follows.
{\hfill \hfill \hfill} \qed

\head 3. Proofs of Theorems 1.3 and 1.6\endhead 

We follow the approach outlined in \cite{Ba16b} and \cite{Ba17}. 
\subhead (3.1) Proof of Theorem 1.3 \endsubhead
Let $A$ be an $n \times n$ complex matrix satisfying the conditions of Theorem  1.3.
We consider a univariate polynomial 
$$g(z) =\per (I + z A).$$
Then $\deg g \leq n$ and, by Theorem 1.2, we have 
$$g(z) \ne 0 \quad \text{provided} \quad |z| < \lambda^{-1}. \tag3.1.1$$
We choose a continuous branch of 
$$f(z) = \ln g(z) \quad \text{for} \quad |z| < \lambda^{-1}. \tag3.1.2$$
Hence our goal is to approximate $f(1)=\ln \per (I + A)$ within an additive error of $\epsilon >0$ by a polynomial $p(A)$ in the entries of $A$ of 
$\deg p =O_{\lambda}\left(\ln n -\ln \epsilon\right)$. 
We consider the Taylor polynomial $T_m(z)$ of $f$ computed at $z=0$:
$$T_m(z)=f(0) + \sum_{k=1}^m {f^{(k)}(0) \over k!} z^k. \tag3.1.3$$
Since $\deg g \leq n$ and (3.1.1) holds, by Lemma 2.2.1 of \cite{Ba16b} (Lemma 7.1 of \cite{Ba17}), we have 
$$\left| f(1) - T_m(1)\right| \ \leq \ {n \lambda^{m+1}  \over (m+1) (1-\lambda)}. \tag3.1.4$$
Therefore, to approximate $f(1)$ within an absolute error $0 < \epsilon < 1$, it suffices to choose $m=O_{\lambda}\left(\ln n - \ln \epsilon\right)$, where the implied constant in the 
``$O$'' notation depends only on $\lambda$.

Next, as is discussed in Section 2.2.2 of \cite{Ba16b} (Section 7.1 of \cite{Ba17}), one can compute the values of $f(0)=\ln g(0)=0$ and $f^{(k)}(0)$ for $k=1, \ldots, m$
from the values of $g(0)=1$ and $g^{(k)}(0)$ for $k=1, \ldots, m$ in $O(m^2)$ time, by solving a non-singular $m \times m$ triangular system of linear equations 
$$\sum_{j=0}^{k-1} {k-1 \choose j} f^{(k-j)}(0) g^{(j)}(0)=g^{(k)}(0) \quad \text{for} \quad k=1, \ldots, m. \tag3.1.5$$
On the other hand, 
$$g(z)=\per(I+ zA) = \sum_{I \subset \{1, \ldots, n\}} z^{|I|} \per A_I$$
where $A_I$ is the submatrix of $A$ consisting of the entries with row and column indices in $I$, where we agree that $\per A_{\emptyset}=1$.
Consequently,
$$g^{(k)}(0)= k! \sum\Sb I \subset \{1, \ldots, n\} \\ |I|=k \endSb \per A_I,$$
so that $g^{(k)}(0)$ is a homogeneous polynomial of degree $k$ in the entries of $A$. It can be computed by the direct enumeration in $n^{O(k)}$ time. It follows from (3.1.5) that $f^{k}(0)$ is a polynomial of degree $k$ in the entries of $A$ computable in $n^{O(k)}$ time. 
Since $m=O_{\lambda}\left( \ln n - \ln \epsilon \right)$, the proof follows by (3.1.4).
{\hfill \hfill \hfill} \qed

\subhead (3.2) Proof of Theorem 1.6 \endsubhead The proof is very similar to that of Theorem 1.3 in Section 3.1. Let $A$ 
be a $d$-dimensional $n \times \ldots \times n$ complex tensor satisfying the conditions of Theorem 1.6. 
We consider the univariate polynomial 
$$g(z) =\PER (I + zA).$$
Then $\deg g \leq n$ and by Theorem 1.5, we have (3.1.1). We then define the function $f(z)$ by (3.1.2), define the polynomial $T_m(z)$ by (3.1.3) and conclude as in Section 3.1 that (3.1.4) holds. As in Section 3.1, to approximate 
$f(1)$ by $T_m(1)$ within an additive error $0 < \epsilon < 1$, it suffices to choose $m=O_{\lambda}(\ln n - \ln \epsilon)$. 

 As in Section 3.1, it remains to prove that $g^{(k)}(0)$ is a polynomial of degree $k$ in the entries of $A$. We have the expansion
 $$\PER (I +zA)=\sum_{I \subset \{1, \ldots, n\}} z^{|I|} \PER A_I,$$
where $A_I$ is the subtensor of $A$ consisting of the entries 
$a_{i_1 \ldots i_d}$ with $i_k \in I$ for $k=1, \ldots, d$ and where we agree that $\PER A_{\emptyset} =1$. Then
$$g^{(k)}(0)=k! \sum\Sb I \subset \{1, \ldots, n\} \\ |I|=k \endSb \PER A_I,$$
 so $g^{(k)}(0)$ is indeed a homogeneous polynomial of degree $k$ in the entries of $A$. Moreover, it can be computed by the direct enumeration in $n^{O(kd)}$ time. Since we have $k=O_{\lambda}(\ln n -\ln \epsilon)$, the proof is completed as in Section 3.1.
{\hfill \hfill \hfill} \qed

\head Acknowledgment \endhead

The author is grateful to Alex Samorodnitsky for suggesting several improvements and to Piyush Srivastava for pointing out to \cite{Br59}.


\Refs
\widestnumber\key{AAAAA}

\ref\key{A+17}
\by N. Anari, L. Gurvits, S. O. Gharan and A. Saberi
\paper Simply exponential approximation of the permanent of positive semidefinite matrices
\paperinfo preprint {\tt arXiv:1704.03486}
\yr 2017
\endref

\ref\key{Ba16a}
\by  A. Barvinok
\paper Computing the permanent of (some) complex matrices
\jour Foundations of Computational Mathematics  
\vol 16 
\yr 2016
\pages  no. 2, 329--342
\endref

\ref\key{Ba16b}
\by  A. Barvinok
\book Combinatorics and Complexity of Partition Functions
\bookinfo Algorithms and Combinatorics, 30
\publ Springer
\publaddr Cham
\yr 2016
\endref

\ref\key{Ba17}
\by A. Barvinok
\paper Approximating permanents and hafnians
\jour Discrete Analysis 
\yr 2017
\pages Paper No. 2, 34 pp
\endref

\ref\key{BR17}
\by A. Barvinok and G. Regts 
\paper Weighted counting of integer points in a subspace
\paperinfo preprint {\tt arXiv:1706.05423}
\yr 2017
\endref

\ref\key{Br59}
\by J.L. Brenner
\paper Relations among the minors of a matrix with dominant principle diagonal
\jour Duke Mathematical Journal 
\vol 26
\yr 1959
\pages 563--567
\endref

\ref\key{EM17}
\by L. Eldar and S. Mehraban
\paper Approximating the permanent of a random matrix with vanishing mean
\paperinfo preprint {\tt arXiv:1711.09457}
\yr 2017
\endref

\ref\key{GS14}
\by  L. Gurvits and A. Samorodnitsky
\paper Bounds on the permanent and some applications
\inbook 55th Annual IEEE Symposium on Foundations of Computer Science -- FOCS 2014
\pages 90--99
\publ IEEE Computer Soc.
\publaddr Los Alamitos, CA
\yr 2014
\endref

\ref\key{J+04}
\by M. Jerrum, A. Sinclair and E. Vigoda
\paper A polynomial-time approximation algorithm for the permanent of a matrix with nonnegative entries
\jour Journal of the ACM 
\vol 51 
\yr 2004
\pages  no. 4, 671--697
\endref

\ref\key{L+17}
\by J. Liu, A. Sinclair and P. Srivastava
\paper The Ising partition function: zeros and deterministic approximation
\paperinfo preprint {\tt arXiv:1704.06493}
\yr 2017
\endref 

\ref\key{MN62}
\by M. Marcus and M. Newman
\paper Inequalities for the permanent function
\jour Annals of Mathematics. Second Series  
\vol 75 
\yr 1962 
\pages 47--62
\endref

\ref\key{PR17a}
\by V. Patel and G. Regts
\paper Deterministic polynomial-time approximation algorithms for partition functions and graph polynomials
\jour  SIAM Journal on Computing
\vol 46 
\yr 2017
\pages no. 6, 1893--1919
\endref

\ref\key{PR17b}
\by H. Peters and G. Regts
\paper On a conjecture of Sokal concerning roots of the independence polynomial
\paperinfo preprint {\tt  arXiv:1701.08049}
\yr 2017
\endref

\endRefs
\enddocument
\end